\documentclass[12pt,twoside]{article}
\usepackage{graphicx}
\usepackage{amsfonts,amsbsy,amssymb,amsmath}
\allowdisplaybreaks
\usepackage{cite}
\usepackage{graphics}
\def\bpsp{\begin{pspicture}}
\def\epsp{\end{pspicture}}

\newtheorem{theorem}{Theorem}[section]
\newtheorem{remark}[theorem]{Remark}
\newtheorem{example}[theorem]{Example}
\newtheorem{lemma}[theorem]{Lemma}
\newtheorem{corollary}[theorem]{Corollary}
\newtheorem{definition}[theorem]{Definition}
\newtheorem{proposition}[theorem]{Proposition}

\newtheorem{note}{Note}
\newtheorem{case}{Case}
\newtheorem{conjecture}{Conjecture}
\newtheorem{question}{Question}

\newcommand{\bea}{\begin{eqnarray}}
\newcommand{\eea}{\end{eqnarray}}
\newcommand{\beq}{\begin{eqnarray*}}
\newcommand{\eeq}{\end{eqnarray*}}

\def\m4{\mbox{\rm ~(mod $4$)}}

\def \bd{\begin{definition}}
\def \ed{\end{definition}}
\def \bqu{\begin{question}}
\def \equ{\end{question}}
\def \bcc{\begin{conjecture}}
\def \ecc{\end{conjecture}}
\def \bt{\begin{theorem}}
\def \et{\end{theorem}}
\def \bl{\begin{lemma}}
\def \el{\end{lemma}}
\def \bc{\begin{corollary}}
\def \ec{\end{corollary}}
\def \be{\begin{equation}}
\def \ee{\end{equation}}
\def \ben{\begin{enumerate}}
\def \een{\end{enumerate}}
\def \ba{\begin{array}}
\def \ea{\end{array}}
\def \bp{\begin{proposition}}
\def \ep{\end{proposition}}
\def \bx{\begin{example}}
\def \ex{\end{example}}
\def \br{\begin{remark}}
\def \er{\end{remark}}
\def \bdsc{\begin{description}}
\def \edsc{\end{description}}

\def \bn{\begin{case}}
\def \en{\end{case}}
\def \bnt{\begin{note}}
\def \ent{\end{note}}
\def\1{1\!\!1}

\def\mm2{\mbox{\rm ~(mod $2$)}}
\def\m4{\mbox{\rm ~(mod $4$)}}

\def\qed{\nolinebreak\hfill\rule{.2cm}{.2cm}\par\addvspace{.5cm}}

\def\m{\mu}

\def\1{\textbf{1}}
\def\0{\textbf{0}}

\begin{document}
\title{On the distribution of eigenvalues of the reciprocal distance Laplacian matrix of graphs}
\author{S. Pirzada$ ^{a} $, Saleem Khan$ ^{b} $ \\
$^{a,b}${\em Department of Mathematics, University of Kashmir, Srinagar, Kashmir, India}\\
$^{a}$\texttt{pirzadasd@kashmiruniversity.ac.in}; $ ^{b} $\texttt{khansaleem1727@gmail.com}
}
\date{}

\pagestyle{myheadings} \markboth{Pirzada,Saleem}{On the distribution of eigenvalues of the reciprocal distance Laplacian matrix of graphs}
\maketitle
\vskip 5mm
\noindent{\footnotesize \bf Abstract.} The Reciprocal distance Laplacian matrix of a connected graph $G$ is defined as $RD^L(G)=RT(G)-RD(G)$,  where $RT(G)$ is the diagonal matrix of reciprocal distance degrees and $RD(G)$ is the Harary matrix. Since $RD^L(G)$ is a real symmetric matrix, we denote its eigenvalues as
  $\lambda_1(RD^L(G))\geq \lambda_2(RD^L(G))\geq \dots \geq \lambda_n(RD^L(G))$. The largest eigenvalue $\lambda_1(RD^L(G))$ of $RD^L(G)$ is called the reciprocal distance Laplacian spectral radius. In this article,  we prove  that the multiplicity of $n$ as a reciprocal distance Laplacian eigenvalue of $RD^L(G)$ is exactly one less than the number of components in the complement graph $\overline{G}$ of $G$. We show that the class of the complete bipartite graphs  maximize the reciprocal distance Laplacian spectral radius among all the bipartite graphs with $n$ vertices. Also, we show that the star graph $S_n$ is the unique graph having the maximum reciprocal distance Laplacian spectral radius in the class of trees with $n$ vertices.  We determine the reciprocal distance Laplacian spectrum of several well known graphs. We prove that the complete graph $K_n$, $K_n-e$, the star $S_n$, the complete balanced bipartite graph $K_{\frac{n}{2},\frac{n}{2}}$ and the complete split graph $CS(n,\alpha)$ are all determined from the $RD^L$-spectrum.

\vskip 3mm



\section{Introduction}

Throughout the paper,  we assume all the graphs under consideration are simple and connected.  A simple connected graph  $G=(V(G),E(G))$ consists of the vertex set $V(G)=\{v_{1},v_{2},\ldots,v_{n}\}$ and the edge set  $E(G)$. The \textit{order} and \textit{size} of $G$ are $|V(G)|=n$ and  $|E(G)|=m$, respectively. The \textit{degree} of a vertex $v_i,$ denoted by $d_{G}(v_i)$ (or shortly $d_i$)  is the number of edges incident on the vertex $v_i$. Further, $N_G (v_i)$ denotes the set of all vertices that are adjacent to $v_i$ in $G$. $\overline{G}$ denotes the complement of the graph $G$.

   The adjacency matrix $A(G)=(a_{ij})$ of $G$ is an $n\times n$ matrix whose $(i,j)$-entry is equal to 1, if $v_i$ is adjacent to $v_j$ and equal to $ 0 $, otherwise. Let $Deg(G)=\text{diag}(d_1, d_2, \dots, d_n)$ be the diagonal matrix of vertex degrees $d_i$, $i=1,2,\dots,n$. The positive semi-definite matrix $L(G)=Deg(G)-A(G)$ is the Laplacian matrix of $G$. The eigenvalues of $L(G)$ are called the Laplacian eigenvalues of $G$. The Laplacian eigenvalues are denoted by $\mu_1 (G) ,\mu_2 (G),\dots,\mu_n (G)$ and are ordered as $\mu_1 (G) \geq \mu_2 (G)\geq \dots \geq \mu_n (G)$. The multiset of  Laplacian eigenvalues of $G$ is called the Laplacian spectrum (briefly $L$-spectrum) of $G$.
   In $G$, the \textit{distance} between two vertices $v_i,v_j\in V(G),$ denoted by $d(v_i,v_j)$, is defined as the length of a shortest path between $v_i$ and $v_j$. The diameter of $G$, denoted by $d(G)$, is $\displaystyle\max_{u,v\in G}d(u,v)$, that is, the length of a longest path among the distance between every two vertices of $G$.  The \textit{distance matrix} of $G$ is denoted by $D(G)$ and is defined as $D(G)=(d(v_i,v_j))_{v_1,v_j\in V(G)}$.
The \textit{transmission} $Tr_{G}(v_i)$
(or briefly $Tr(i)$ if graph $G$ is understood) of a vertex $v_i$ is defined as the sum of the distances from $v_i$ to all other vertices in $G$, that is, $Tr_{G}(v_i)=\sum\limits_{v_j\in V(G)}d(v_i.v_j).$ For any vertex $v_i\in V(G)$, the transmission $Tr_G(v_i)$ is also called the \textit{transmission degree}.

Let $Tr(G)=diag (Tr_1,Tr_2,\ldots,Tr_n) $ be the diagonal matrix of vertex transmissions of $G$. Aouchiche and Hansen \cite{9R6} introduced the Laplacian for the distance matrix of a connected graph. The matrix $D^L(G)=Tr(G)-D(G)$ (or simply $D^{L}$) is called the \textit{distance Laplacian matrix} of $G$. The eigenvalues of $D^{L}(G)$ are called the distance Laplacian eigenvalues of $G$ and are referred as $D^L$-eigenvalues of $G$. Since $ D^L(G) $ is a real symmetric positive semi-definite matrix, we take the distance Laplacian eigenvalues as $\partial_{1}^{L}(G)\geq \partial_{2}^{L}(G)\geq \dots\geq \partial_{n}^{L}(G)$.

The \textit{Harary matrix} of graph $G$, which is also called as the Reciprocal Distance matrix, denoted by $RD(G)$, is an $n$ by $n$ matrix defined as \cite{9R8}
\begin{equation*}
RD_{ij}=\begin{cases}
 \frac{1}{d(v_i,v_j)}& \text{if} ~ i\neq j\\
 0 &\text{if} ~ i=j.\\
\end{cases}
\end{equation*}

Henceforward, we consider $i\neq j$ for $d(v_i,v_j)$.

The Reciprocal distance degree of a vertex $v_i$, denoted by $RTr_G(v_i)$ (or shortly $RTr(i)$ ), is given by
$$RTr_G(v_i)=\displaystyle\sum_{v_j\in V(G)_{v_i\neq v_j}}\frac{1}{d(v_i,v_j)}.$$
Let $RT(G)$ be an $n\times n$ diagonal matrix defined by $RT_{ii}=RTr_G(v_i)$.

The Harary index of a graph $G$, denoted by $H(G)$, is defined in \cite{9R8} as

$$H(G)=\frac{1}{2}\displaystyle\sum_{i=1}^{n}\displaystyle\sum_{j=1}^{n}RD_{ij}=\frac{1}{2}\displaystyle\sum_{v_j\in V(G)_{v_i\neq v_j}}\frac{1}{d(v_i,v_j)}.$$
Clearly,
$$H(G)=\frac{1}{2}\displaystyle\sum_{v_i\in V(G)}RTr_G(v_i).$$
 To see more work done on the Harary matrix, we refer to \cite{9R9,9R10,9R11}.

 In \cite{9R1}, the authors defined the Reciprocal distance Laplacian matrix as $RD^L(G)=RT(G)-RD(G)$. Since $RD(G)$ and $RD^L(G)$ are real symmetric matrices, we can denote by
 $$\lambda_1(RD(G))\geq \lambda_2(RD(G))\geq \dots \geq \lambda_n(RD(G)),$$
 and
  $$\lambda_1(RD^L(G))\geq \lambda_2(RD^L(G))\geq \dots \geq \lambda_n(RD^L(G))$$
the eigenvalues of $RD(G)$ and $RD^L(G)$, respectively. Since $RL(G)$ is a positive semidefinite matrix ,we will denote by $\rho(RD^L(G))=\lambda_1(RD^L(G))$ the spectral radius of $RL(G)$. In \cite{9R1}, the authors proved that for the connected graph $G$ of order $n$, the spectral radius of $RD^L(G)$  is at most $n$ and supplied the necessary and sufficient conditions for $n$ to be the eigenvalue of  $RD^L(G)$. More work on the matrix  $RD^L(G)$ can be seen in \cite{9R12, 9R13, 9R14}.

The rest of the paper is organised as follows. In Section 2, we supply some results from the previous works and  prove that the multiplicity of $n$ as a reciprocal distance Laplacian eigenvalue of $RD^L(G)$ is exactly one less than the number of components in the complement graph $\overline{G}$ of $G$. We show that the class of the complete bipartite graphs  maximize the reciprocal distance Laplacian spectral radius among all the bipartite graphs. We also show that the star graph $S_n$ is the unique graph having the maximum reciprocal distance laplacian spectral radius in the class of trees. In Section 3, we find the reciprocal distance Laplacian spectrum of several well known graphs. In Section 4, we prove that the complete graph $K_n$, $K_n-e$, the star $S_n$, the complete balanced bipartite graph $K_{\frac{n}{2},\frac{n}{2}}$ and the complete split graph $CS(n,\alpha)$ are all determined from the $RD^L$-spectrum.

\section{On the eigenvalues of reciprocal distance Laplacian matrix of graph $G$}

We begin this section with the following lemmas which will be used in sequel.
\begin{lemma}\label{L1}\emph{\cite{9R2}}
If the graph $G$ has $n$ vertices and $\mu$ is an eigenvalue of $L(G)$, then $0\leq \mu \leq n$. Moreover, the multiplicity of $n$ is equal to one less than the number of components in the complement graph $\overline{G}$.
\end{lemma}
\begin{lemma}\label{L2}\emph{\cite{9R1}}
Let $G$ be a connected graph on $n$ vertices with diameter $d=2$. Then $\lambda_i(RD^L(G))=\frac{n+\mu_i(G)}{2}$ for $i=1,2,\dots,n-1$. Furthermore, $\frac{n+\mu_i(G)}{2}$ and $\mu_i(G)$ both have the same multiplicity for $i=1,2,\dots,n$.
\end{lemma}
\begin{lemma}\label{L3}\emph{\cite{9R1}}
Let $G$ be a connected graph on $n$ vertices. Then, the complement graph $\overline{G}$ is disconnected if and only if the reciprocal distance Laplacian spectral radius of $G$ is $n$.
\end{lemma}
\begin{lemma}\label{L11}\emph{\cite{9R1}}
For any connected graph $G$, $0$ is a simple eigenvalue of $RD^L(G)$.
\end{lemma}
\begin{lemma}\label{L4}\emph{\cite{9R1}}
Let $G$ be a connected graph and $G'=G+e$, where $e\not\in E(G)$. Then $\lambda_i(RD^L(G'))\geq \lambda_i(RD^L(G))$ for all $i=1,2,\dots,n$.
\end{lemma}
The next lemma follows immediately from the above lemma.
\begin{lemma}\label{L5}
Let $G$ be a  connected graph on $n$ vertices. Then
$$\lambda_n(RD^L(G))=0~ \text{and}~ \lambda_i(RD^L(G))\leq \lambda_i(RD^L(K_n))=n ~\text{for  all}~ i=1,2,\dots,n-1.$$
\end{lemma}
In the following theorem, we show that $n$ as an eigenvalue of $RD^L(G)$ can be considered as the algebraic connectivity of $\overline{G}$, that is, we show that the multiplicity of $n$  as an eigenvalue of $RD^L(G)$ is exactly one less than the number of components in the complement graph $\overline{G}$. We note here that the same type of result has been proved for the distance Laplacian matrix in \cite{9R6}
\begin{theorem}\label{T1}
For every connected graph on $n$ vertices, $n$ is an eigenvalue of $RD^L(G)$ with multiplicity exactly  equal to one less than the number of components in the complement graph $\overline{G}$.
\end{theorem}
\noindent{\bf Proof.} If the complement graph $\overline{G}$ of $G$ is connected, then by Lemma \ref{L3}, the result is trivially true. So, let the complement graph $\overline{G}$ be disconnected so that diameter of $G$ is 2. By Lemma \ref{L2}, the eigenvalues of $RD^L(G)$ are given by $\lambda_i(RD^L(G))=\frac{n+\mu_i(G)}{2}$ for $i=1,2,\dots,n-1$ and $\lambda_n(RD^L(G))=0$, where $\mu_1 \geq \mu_2 \geq \dots \geq \mu_{n-1}$ are the Laplacian eigenvalues of $G$. Using Lemma \ref{L1}, we get the desired result. \qed
The following observations are an immediate consequence of Theorem \ref{T1}.
\begin{corollary}\label{C1}
Let $G$ be a connected graph on $n$ vertices. Then, $\lambda_{n-1}(RD^L(G))\leq n$ with equality if and only if $G\cong K_n$.
\end{corollary}
\noindent{\bf Proof.}  The bound follows by Lemma \ref{L5}. Further, we observe that equality holds for $K_n$.

Let $\lambda_{n-1}(RD^L(G))=n$. By Lemma \ref{L5}, $n$ is an eigenvalue of $RD^L(G)$ with multiplicity $n-1$. Thus , by Theorem \ref{T1}, $\overline{G}$ has $n$ components which are necessarily isolated vertices and therefore $G$ is a complete graph. \qed
\begin{corollary}\label{C2}
Let $G\ncong K_n$ be a connected graph on $n$ vertices and let $D=\{v\in V(G):d_G( v)=n-1\}$. Then $n$ is an eigenvalue of $RD^L(G)$ with multiplicity at least $|D|$.
\end{corollary}
\noindent{\bf Proof.} We note that each vertex $v\in D$ of $G$ corresponds to an isolated vertex in $\overline{G}$ and thus to a component of $\overline{G}$. As $G\ncong K_n$, the number of components in the complement graph $\overline{G}$ is at least $|D|+1$ and the result follows by using Theorem \ref{T1}. \qed

Let $\mathbb{B}_n $ be the set of all connected bipartite graphs on $n$ vertices and $\mathbb{K}^b_a$ be the set of all complete bipartite graphs on $n=a+b$ vertices. Let $\mathbb{T}_n$ be the set of all trees on $n$ vertices. The following result characterizes all the graphs in $\mathbb{B}_n $ and $\mathbb{T}_n$ having the maximum reciprocal distance Laplacian spectral radius.

\begin{theorem}\label{T2} Let $G$ be a connected graph on $n$ vertices.

 (a) If $G\in \mathbb{B}_n $, then $\rho(RD^L(G))\leq \rho(RD^L(K_{a,b}))$ , where $K_{a,b}$ is any graph in $\mathbb{K}^b_a$,  with equality if and only if $G\in \mathbb{K}^b_a$.

(b) If $G\in \mathbb{T}_n$, then $\rho(RD^L(G))\leq  \rho(RD^L(S_n)) $ with equality if and only if $G\cong S_n$.

\end{theorem}
\noindent{Proof.} (a) We observe that the complete bipartite graphs are the only connected bipartite graphs that have disconnected complement. Thus, from Theorem \ref{T1}, the complete bipartite graphs are the only connected bipartite graphs that have $n$ as a reciprocal  distance Laplacian eigenvalue and the result follows from Lemma \ref{L5}.

\noindent (b) We see that the star $S_n$ is the only tree on $n$ vertices having disconnected complement which by Theorem \ref{T1} shows that among all the trees on $n$ vertices only the star $S_n$ has $n$  as a reciprocal  distance Laplacian eigenvalue and the result follows from Lemma \ref{L5}.\qed

\section{Reciprocal Distance Laplacian spectrum}

It is possible to know some reciprocal distance Laplacian eigenvalues of a graph $G$ given that $G$ has some particular structure.
 \begin{theorem}\label{T3}
Let $G$ be a connected graph on $n\geq 2$ vertices. If $M=\{v_1,v_2,\dots,v_r\}$ is an independent set of $G$ such that $N_G(v_i)= N_G(v_j)$ for all $1\leq i,j\leq r$, then $T=RTr(v_i)=RTr(v_j)$  for all $1\leq i,j\leq r$ and $T+\frac{1}{2}$ is an eigenvalue of $RD^L(G)$ with multiplicity at least $r-1$.
 \end{theorem}
\noindent{\bf Proof.} We see that any two vertices in $M$ are at a distance of 2 from each other. So $RD_{ij}=\frac{1}{2}$, for all $1\leq i\neq j\leq r$. Also any vertex in $V(G)\setminus M$ is at the same distance from all the vertices in $M$. Thus all the vertices in $M$ have the same reciprocal distance degree, say $T$.

 The proof gets completed after observing that the matrix $(T+\frac{1}{2})I_n-RD^L(G)$ has at least $r$ identical rows. \qed

By using the similar argument, we can prove the following theorem.

\begin{theorem}\label{T4}
Let $G$ be a connected graph on $n\geq 2$ vertices. If $S=\{v_1,v_2,\dots,v_r\}$, is a clique of $G$ such that $N_G(v_i)-S= N_G(v_j)-S$ for all $1\leq i,j\leq r$, then $T=RTr(v_i)=RTr(v_j)$  for all $1\leq i,j\leq r$ and $T+1$ is an eigenvalue of $RD^L(G)$ with multiplicity at least $r-1$.
 \end{theorem}

We note that the results similar to Theorems \ref{T3} and \ref{T4} have been proved for the distance Laplacian matrix (see \cite{9R4}) and  distance signless Laplacian matrix (see \cite{9R3}).

By using Theorems \ref{T3} and \ref{T4}, we obtain the reciprocal distance Laplacian spectrum of some well known graphs in the following result.
 \begin{lemma}\label{L6}
 (a)~  The  reciprocal distance Laplacian characteristic polynomial of the complete $k$-partite graph $K_{n_1,n_2,\dots,n_k}$, where $n_1+n_2+\dots+n_k=n$, is
 $$P^x_{RD^L}(K_{n_1,n_2,\dots,n_k})=x(x-n)^{k-1}\displaystyle\prod_{i=1}^{k}\Bigg(x-\Big(n-\frac{n_i}{2}\Big)\Bigg)^{n_i-1}.$$
  (b)~  The  reciprocal distance Laplacian characteristic polynomial of the complete bipartite graph $K_{a,b}$, where $a+b=n$, is
   $$P^x_{RD^L}(K_{a,b})=x(x-n)\Bigg(x-\Big(n-\frac{a}{2}\Big)\Bigg)^{a-1}\Bigg(x-\Big(n-\frac{b}{2}\Big)\Bigg)^{b-1}.$$
    (c)~  The  reciprocal distance Laplacian characteristic polynomial of the star $S_n$ is
    $$P^x_{RD^L}(S_n)=x(x-n)\Big(x-\frac{n+1}{2}\Big)^{n-2}.$$
    (d)~  The  reciprocal distance Laplacian characteristic polynomial of the complete split graph $CS(n,\alpha)$ is given by
    $$P^x_{RD^L}(CS(n,\alpha))=x(x-n)^{n-\alpha}\Bigg(x-\Big(n-\frac{\alpha}{2}\Big)\Bigg)^{\alpha-1}.$$
    (e)~  The  reciprocal distance Laplacian characteristic polynomial of the graph $K_n-e$ is
     $$P^x_{RD^L}(K_n-e)=x(x-n+1)(x-n)^{n-2}.$$
    (f)~  The  reciprocal distance Laplacian characteristic polynomial of the graph $PA(n,p)$, obtained from a clique $K_{n-p}$ by attaching $p>0$ pendant edges to a vertex from the clique, is
     $$P^x_{RD^L}(PA(n,p))=x(x-n)\Big(x-n+\frac{p}{2}\Big)^{n-p-2}\Big(x-\frac{n+1}{2}\Big)^p.$$
   (g)~  The  reciprocal distance Laplacian characteristic polynomial of the graph $S^+_n$, obtained from the star $S_n$ by adding an edge,  is
     $$P^x_{RD^L}(S^+_n)= x(x-n)\Big(x-\frac{n+3}{2}\Big)\Big(x-\frac{n+1}{2}\Big)^{n-3}.$$
 \end{lemma}

\noindent{\bf Proof.} (a)~ The complete $k$-partite graph $K_{n_1,n_2,\dots,n_k}$ contains $k$ independent sets and let $M_i$, $i=1,2,\dots,k$, be these $k$ independent sets with respective cardinalities $|M_i |=n_i$, $i=1,2,\dots,k$. For each $i=1,2,\dots,k$,  every vertex in $M_i$ share the same neighbourhood $V(K_{n_1,n_2,\dots,n_k})\setminus M_i$ and each vertex in $M_i$ is at a distance of 2 from every vertex in $M_i$. Thus, the reciprocal distance degree of each vertex in $M_i$ is $\frac{2n-n_i-1}{2}$, $i=1,2,\dots,k$. By Theorem \ref{T3}, $\frac{2n-n_i}{2}$ is an eigenvalue of $RD^L(G)$ with multiplicity at least $n_i -1$, $i=1,2,\dots,k$. In addition, the complement of the complete $k$-partite graph $K_{n_1,n_2,\dots,n_k}$ contains exactly $k$ components so that $n$ is an  eigenvalue of $RD^L(G)$ with multiplicity exactly $k-1$ and the remaining one eigenvalue by Lemma \ref{L11} is 0.

\noindent (b)~  Putting $k=2$, $n_1=a$ and $n_2=b$ in the complete $k$-partite graph $K_{n_1,n_2,\dots,n_k}$, we get the required result.

\noindent (c)~ Putting $a=1$ and $b=n-1$ in (b) gives the desired result.

\noindent (d)~ Clearly, $CS(n,\alpha)$ has $n-\alpha+1$ independent sets. The largest one has cardinality $\alpha$ and the remaining are of cardinality 1. So, putting $k=n-\alpha+1$, $n_1=\alpha$ and  $n_2=n_3=\dots=n_{n-\alpha+1}=1$ in (a) proves the result.

\noindent (e)~ We observe that the complement of $K_n-e$ contains $n-1$ components out of which $n-2$ components are isolated vertices and one component is $K_2$. By Theorem \ref{T1}, $n$ is an eigenvalue of $RD^L(K_n-e)$ with multiplicity $n-2$. Using Theorem \ref{T3}, the remaining non-zero eigenvalue of $RD^L(K_n-e)$  is seen to be $n-1$.

\noindent (f)~ Clearly, $PA(n,p)$ contains a clique on $n-p-1$ sharing the same neighbourhood (a dominating vertex) and the same reciprocal distance degree $n-\frac{p}{2}-1$. By Theorem \ref{T4}, $n-\frac{p}{2}$ is an eigenvalue of $RD^L(PA(n,p))$ with multiplicity at least $n-p-2$.  Als0, $PA(n,p)$ contains an independent set of $p$  vertices  sharing the same neighbourhood and the same reciprocal distance degree $\frac{n}{2}$. Thus, by Theorem \ref{T3}, $\frac{n+1}{2}$ is an eigenvalue of  $RD^L(PA(n,p))$ with  multiplicity at least $p-1$. Now, using the facts that 0 is always a simple eigenvalue and the sum of all reciprocal distance degrees is equal to sum of all reciprocal distance Laplacian eigenvalues, we can easily evaluate the remaining eigenvalue, which equals $\frac{n+1}{2}$ .

\noindent (g)~ Putting $p=n-3$ in (f), we get the required result.\qed

   We also make use of the following lemmas.

  \begin{lemma}\label{L7} \emph{\cite{9R2}}
  Let $C_n$ and $P_n$ be the cycle and path on $n$ vertices, respectively. The Laplacian eigenvalues of $C_n$ and $P_n$ are given by $4sin^2\frac{\pi k}{n}$, $k=1,2,\dots,n$ and $4sin^2\frac{\pi k}{2n}$, $k=0,1,2,\dots,n-1$, respectively.
  \end{lemma}
  \begin{lemma}\label{L12}\emph{\cite{9R5}}
  Let $G$  be a graph on $n$ vertices. If $\mu_i(G)$, $i=1,2,\dots,n$ are the eigenvalues of $L(G)$ then the eigenvalues of $L(\overline{G})$ are $n-\mu_{n-i}(G)$, $i=1,2,\dots,n-1$ and $0$.
  \end{lemma}

  The next observation follows immediately from Lemmas \ref{L7} and \ref{L12}. This gives the reciprocal distance Laplacian eigenvalues of the complements of  path and cycle.

  \begin{lemma}\label{L8}  For $n\geq 5$,
    the reciprocal distance Laplacian characteristic polynomial of $\overline{C}_n$ is given by
$$x\displaystyle\prod_{i=1}^{n-1}\Bigg(x-n+2sin^2\Big(\frac{\pi i}{n}\Big)\Bigg)$$
and
  the reciprocal distance Laplacian characteristic polynomial of $\overline{P}_n$ is given by
$$x\displaystyle\prod_{i=1}^{n-1}\Bigg(x-n+2sin^2\Big(\frac{\pi i}{2n}\Big)\Bigg).$$
  \end{lemma}
\noindent{\bf Proof.} Using Lemmas \ref{L7} and \ref{L12}, we see that the Laplacian eigenvalues of $\overline{C}_n$  are given by  $n-4sin^2\frac{\pi k}{n}$, $k=1,2,\dots,n-1$ and $0$. Similarly, the Laplacian eigenvalues of $\overline{P}_n$  are given by $n-4sin^2\frac{\pi k}{2n}$, $k=1,2,\dots,n-1$ and $0$. For $n\geq 5$, it is easy to see that both $\overline{C}_n$ and $\overline{P}_n$ are  connected with diameter 2 and the proof follows from Lemma \ref{L2}. \qed

\section{Graphs determined by $RD^L$- spectrum}

 In this section, we show the existence of some graphs which are determined by $RD^L$- spectrum. We start with the following observations.

 \begin{lemma}\label{L9}\emph{\cite{9R6}} Let $G$ be a connected graph on $n$ vertices with $diam(G)\leq 2$. Let $\mu_1 (G) \geq \mu_2 (G)\geq \dots \geq \mu_n (G)=0$ be the Laplacian spectrum of $G$. Then the distance Laplacian spectrum of $G$ is  $2n-\mu_{n-1} (G) \geq 2n- \mu_{n-2} (G)\geq \dots \geq 2n-\mu_1 (G)>\partial^L_n (G)=0$. Moreover, for every $i\in \{1,2,\dots,n-1\}$ the eigenspaces corresponding to $\mu_i (G)$ and $2n-\mu_i (G)$ are same.
\end{lemma}
 \begin{lemma}\label{L10}\emph{\cite{9R7}}
  Let $G$ be a connected graph on $n$ vertices such that $D^L(G)$ has an eigenvalue with multiplicity $n-2$. Then, $m(\partial^L_1(G))=n-2$ if and only if $G\cong S_n$ or $G\cong K_{\frac{n}{2},\frac{n}{2}}$, if $n$ is even.
 \end{lemma}
 \begin{theorem}\label{T5}
 From the reciprocal distance Laplacian spectrum of the connected graph $G$, we can determine (a)~ The number of vertices of $G$, (b)~ the Harary index of $G$ and (c)~ the number of components of complement graph $\overline{G}$.
 \end{theorem}
\noindent{\bf Proof.} (a)~ The number of vertices of $G$ is equal to the number of eigenvalues.

\noindent (b)~ The Harary index of $G$ is half the sum of the reciprocal distance degrees of the vertices of $G$, which is half of the sum of the eigenvalues of $RD^L(G)$.

\noindent (c)~ By Theorem \ref{T1}, the multiplicity  of $n$ as an eigenvalue of the matrix $RD^L(G)$ is exactly one less than the number of components in the complement graph $\overline{G}$.\qed

 If $G$ is a connected graph with diameter 2, then the above theorem  can be improved as follows.
  \begin{theorem}\label{T6}
 From the reciprocal distance Laplacian spectrum of the connected graph $G$ having diameter 2, we can determine (a)~ the number of the vertices of $G$, (b)~ the Harary index of $G$, (c)~ the number of the components of the complement graph $\overline{G}$, (d)~ the Laplacian eigenvalues of $G$ including their multiplicities and (e)~ the number of the edges in $G$.
 \end{theorem}
\noindent{\bf Proof.} (a), (b) and (c) are proved in Theorem \ref{T5}.

\noindent (d)~ It follows directly from Theorem \ref{T1}.

\noindent (e)~ The number of edges of  $G$ is half the sum of the vertex degrees of $G$, which is half of the sum of the Laplacian  eigenvalues of $ G$ and the proof follows from (d). \qed

 \begin{lemma} \label{L12}
 Let $G$ be a connected graph on $n$ vertices and let $H(G)$ be the Harary index of $G$. Then,

 (a)~ $H(G)\leq \frac{n(n-1)}{2}$ with equality if and only if $G\cong K_n$;

 (b)~ $H(G)\leq \frac{n(n-1)}{2}-\frac{1}{•2}$ with equality if and only if $G\cong K_n- e$.
 \end{lemma}
\noindent{\bf Proof.} (a)~ We know that $H(G)=\frac{1}{2}\displaystyle\sum_{u,v\in V(G)_{u\neq v}}\frac{1}{d(u,v)}$. Since the distance between any two different vertices in $K_n$ is equal to 1, therefore, $$H(K_n)=\frac{1}{2}\displaystyle\sum_{u,v\in V(G)_{u\neq v}}\frac{1}{d(u,v)}=\frac{1}{2}\displaystyle\sum_{u,v\in V(G)_{u\neq v}}1=\frac{n(n-1)}{2}.$$

  Thus, the equality holds for $K_n$.
  Now, suppose that $	G\ncong K_n$. Then there are at least two vertices which are non-adjacent, say $u$ and $v$, and $d(u,v) \geq 2$ so that $\frac{1}{d(u,v)}\leq \frac{1}{2} $. The rest of the distances in $G$ are at least 1 so that their reciprocal is at most 1. Therefore,
  $$H(G)=\frac{1}{2}\displaystyle\sum_{u,v\in V(G)_{u\neq v}}\frac{1}{d(u,v)}\leq \frac{1}{2}\Bigg(n(n-1)-2+2\Big(\frac{1}{2}\Big)\Bigg)=\frac{n(n-1)}{2}-\frac{1}{2}.$$

\noindent  (b)~ The graph $K_n -e$ is obtained uniquely from $K_n$ by deleting a single edge from  $K_n$ and this operation decreases the value of   $H(K_n)$ by $\frac{1}{2}$. Removal of any further edges strictly decreases the value of the Harary index. Thus, $K_n -e$ is the unique graph having Harary index equal to $\frac{n(n-1)}{2}-\frac{1}{•2}$. \qed

 \begin{theorem}\label{T7}  The following graphs are determined by their reciprocal distance Laplacian spectrum.

 (a)~ the complete graph $K_n$,

 (b)~ the graph $K_n-e$ obtained from $K_n$ by the deletion of an edge,

 (c)~ the star $S_n$,

 (d)~ the complete balanced bipartite graph $K_{\frac{n}{2},\frac{n}{2}}$,

 (e)~ the complete split graph $CS(n,\alpha)$.
  \end{theorem}

\noindent{\bf Proof.} (a)~ Let $G$ be the graph with the same reciprocal distance Laplacian spectrum as $K_n$, which is given by $\{0,n^{(n-1)}\}$. Since $n$ is $RD^L(G)$-eigenvalue of $G$ with multiplicity $n-1$, therefore, by Theorem \ref{T1}, the complement graph $\overline{G}$ of $G$ has $n$ components that necessarily are isolated vertices which is only possible if $G$ is the complete graph $K_n$.\\
\noindent (b)~ It follows directly from Theorem \ref{T5} and Lemma \ref{L12}.\\
\noindent (c)~ Let $G$ be the graph with the same reciprocal distance Laplacian spectrum as $S_n$, which from Lemma \ref{L6} is given by $\Big\{0,n, {\frac{n+1}{2}}^{(n-2)}\Big\}$. As $n$ is an eigenvalue of $RD^L(G)$ with multiplicity one, therefore, by Theorem \ref{T1}, the complement graph $\overline{G}$ of $G$ has two components. This shows that the diameter of $G$ is 2. Thus, from Theorem \ref{T6}, we can find the number of edges in $G$ which is $n-1$. This shows that $G$ must be a tree on $n$ vertices. Combining the facts that $G$ is a tree and its complement is disconnected, we observe that $G$ is isomorphic to $S_n$, as $S_n$ is the only tree whose complement is disconnected.\\
\noindent (d)~ Let $G$ be the graph with the same reciprocal distance Laplacian spectrum as $K_{\frac{n}{2},\frac{n}{2}}$, which from Lemma \ref{L6} is given by $\Big\{0,n, {\frac{3n}{4}}^{(n-2)}\Big\}$.  As $n$ is an eigenvalue of $RD^L(G)$ with multiplicity one, therefore, by Theorem \ref{T1}, the complement graph $\overline{G}$ of $G$ has two components which shows that diameter of $G$ is 2. Thus, from Theorem \ref{T6}, we can determine the  Laplacian eigenvalues of  $G$ which are given by $\Big\{0,n, {\frac{n}{2}}^{(n-2)}\Big\}$ and the number of edges in $G$ which is $\frac{n^2}{4}$. From Lemma \ref{L9}, the distance Laplacian spectrum of $G$ is given by  $\Big\{0,n, {\frac{3 n}{2}}^{(n-2)}\Big\}$, so that $m(\partial^L_1(G))=n-2$. Using Lemma \ref{L10}, $G\cong S_n$ or $G\cong K_{\frac{n}{2},\frac{n}{2}}$. Using the facts that $G$ has $\frac{n^2}{4}$ edges and $S_n$ has only $n-1$ edges , we see that $G\cong K_{\frac{n}{2},\frac{n}{2}}$.\\
\noindent (e)~ Let $G$ be the graph with the same reciprocal distance Laplacian spectrum as $CS(n,\alpha)$ which from Lemma \ref{L6} is given by $\Big\{0,n^{(n-\alpha)}, {n-\frac{\alpha}{2}}^{(\alpha-1)}\Big\}$.  As $n$ is an eigenvalue of $RD^L(G)$ with multiplicity $n-\alpha$, therefore, by Theorem \ref{T1}, the complement graph $\overline{G}$ of $G$ has exactly $n-\alpha+1$ components. This shows that the diameter of $G$ is 2. Using  Lemma \ref{L2}, the Laplacian eigenvalues of ${G}$ are given by $\{0,n^{(n-\alpha)},{n-\alpha}^{(\alpha-1)}\}$. Using Lemma \ref{L12}, the Laplacian eigenvalues of the complement graph $\overline{G}$  of $G$ are given by $\{0^{(n-\alpha+1)}, {\alpha}^{(\alpha-1)}\}$. As $\overline{G}$ has exactly $n-\alpha+1$ components, therefore, the largest component in $\overline{G}$ contains at most $\alpha$ vertices. We claim that at least one  component in $\overline{G}$ contains exactly $\alpha $ vertices. If possible, let all the  components in $\overline{G}$, say $H_i$ ~ $(1\leq i \leq n-\alpha+1)$, contain less than $\alpha$ vertices. Using the fact that the Laplacian spectral radius of any connected graph  is always less than or equal to the order of the graph, we get $\mu_1(H_i)< \alpha$ ~ $(1\leq i \leq n-\alpha+1)$. This shows that $\mu_1(\overline{G})<\alpha$. This is a contradiction, since from the Laplacian spectrum of $\overline{G}$, we have $\mu_1(\overline{G})=\alpha$. That proves the claim. Since the complement graph $\overline{G}$ of $G$ has exactly $n-\alpha+1$ components, we see that exactly one component of   $\overline{G}$ contains $\alpha$ vertices and the rest of the components are isolated vertices. Since multiplicity of $\alpha$ as a Laplacian eigenvalue of $\overline{G}$ is $\alpha-1$, therefore, $\overline{G}$ must be isomorphic to $K_{\alpha}\cup (n-\alpha)K_1$, which further shows that $G\cong CS(n,\alpha)$. \qed

\noindent{\bf Acknowledgements.} The research of S. Pirzada is supported by SERB-DST, New Delhi under the research project number CRG/2020/000109.\\

\noindent{\bf Data availibility} Data sharing is not applicable to this article as no data sets were generated or analyzed
during the current study.\\	

\noindent{\bf Conflict of interest.} The authors declare that they have no conflict of interest.


\begin{thebibliography}{99}

\bibitem{9R6} M. Aouchiche and P. Hansen, Two Laplacians for the distance matrix of a graph, Linear Algebra Appl. {\bf 439 (1)} (2013) 21-33.
\bibitem{9R4} M. Aouchiche and P. Hansen, Some properties of the distance Laplacian eigenvalues of a graph, Czechoslovak Mathematical Journal {\bf 64 (139)} (2014) 751-761.
\bibitem{9R3} M. Aouchiche and P. Hansen, On the distance signless Laplacian of a graph, Linear Multilinear Algebra {\bf 64} (2016) 1113-1123.
\bibitem{9R2} W. N. Anderson and T. D. Morley, Eigenvalues of the Laplacian of a graph, Linear Multilinear Algebra {\bf 18(2)} (1985) 141-145.
\bibitem{9R1} R. Bapat and S. K. Panda, The spectral radius of the reciprocal distance Laplacian matrix of a graph, Bull. Iranian Math. Soc. (2018) {\bf 44(5)} 1211-1216.
\bibitem{9R9} K. C. Das, Maximum eigenvalue of the reciprocal distance matrix, J Math. Chem. {\bf 47} (2010) 21-28.
\bibitem{9R7} R. Fernandes, M. Aguieiras, A. de Freitas, C. M. da Silva Jr. anf R. R. Del-Vecchio, Multiplicities of distance Laplacian eigenvalues and forbidden subgraphs, Linear Algebra Appl. {\bf 541} (2018) 81-93.
\bibitem{9R11} F. Huang, X. Li and S. Wang, On graphs with maximum Harary spectral radius, (2014). Available from arXiv:1411.6832v1.
\bibitem{9R12} L. Medina and M. Trigo, Upper bounds and lower bounds for the spectral radius of reciprocal distance, reciprocal distance Laplacian and reciprocal distance signless Laplacian matrices, Linear Algebra Appl. {\bf 609} (2021) 386-412.
\bibitem{9R13} L. Medina and M. Trigo, Bounds on the  reciprocal distance energy and  reciprocal distance Laplacian energies of a graph, Linear Multilinear Algebra, DOI: 10.1080/03081087.2020.1825607.
\bibitem{9R5} R. Merris, Laplacian matrices of graphs: A survey, Linear Algebra Appl. {\bf 197-198} (1994) 143-176.
\bibitem{9R8} D. Plavsi$\acute{c}$, S. Nikoli$\acute{c}$, N. Trinajsti$\acute{c}$ and Zlatko Mihali$\acute{c}$, On the Harary index for the characterization of chemical graphs, J. Math. Chem. {\bf 12} (1993) 235-250.
\bibitem{9R14} M. Trigo, On Hararay energy and reciprocal distance Laplacian energies, J. Phys.: Conf. Ser. {\bf 2090} (2021) 012102.
\bibitem{9R10} B. Zhou,  N. Trinajsti$\acute{c}$, Maximum eigenvalues of the reciprocal distance matrix and the reverse Wiener matrix. Int. J. Quantum Chem. \textbf{108} (2008) 858-864.
\end{thebibliography}
\end{document}